\documentclass[11pt]{article}
\usepackage{amsmath,amssymb,amsfonts,amsthm,mathtools}
\usepackage[top=25mm,bottom=25mm,left=25mm,right=25mm]{geometry}
\usepackage{microtype,booktabs}
\usepackage{tikz}
\usepackage{comment}
\usetikzlibrary{arrows.meta}
\usepackage[colorlinks=true,linkcolor=blue,citecolor=blue,urlcolor=cyan]{hyperref}
\usepackage[nameinlink]{cleveref}
\theoremstyle{plain}
\newtheorem{theorem}{Theorem}[section]
\newtheorem{lemma}[theorem]{Lemma}
\newtheorem{proposition}[theorem]{Proposition}
\newtheorem{corollary}[theorem]{Corollary}
\theoremstyle{definition}
\newtheorem{example}[theorem]{Example}
\newtheorem{remark}[theorem]{Remark}
\newtheorem{problem}[theorem]{Problem}
\newtheorem{innerthm}{Theorem}
\newenvironment{reptheorem}[1]{\renewcommand\theinnerthm{\ref*{#1}}\innerthm}{\endinnerthm}
\newtheorem{innerprop}{Proposition}
\newenvironment{repproposition}[1]{\renewcommand\theinnerprop{\ref*{#1}}\innerprop}{\endinnerprop}
\newcommand{\F}{\mathbb F}
\newcommand{\abs}[1]{\lvert#1\rvert}
\newcommand{\set}[1]{\{#1\}}
\newcommand{\angles}[1]{\langle#1\rangle}
\DeclareMathOperator{\Stab}{Stab}
\title{The additive triple spectrum in $\F_3^n$}
\author{
Laura Johnson\thanks{School of Mathematics, University of Bristol, England, UK \emph{E-mail}: \text{laura.marie.johnson@bristol.ac.uk}.} 
\and 
Ananthakrishnan Ravi\thanks{Delft University of Technology, Netherlands. Supported by an NWO open competition
grant (OCENW.M.22.090). \emph{E-mail}: \text{a.ravi@tudelft.nl}.}
}
\date{}
\begin{document}
\maketitle

\begin{abstract}
For subsets $A,B\subseteq\mathbb F_3^n$ of sizes $s,t$, respectively,
let $r(A,B,B)$ be the number of ordered triples
$(a,b,b')\in A \times B\times B$ satisfying $a+b=b'$.
The additive triple spectrum $r_n(s,t)$ is the set of all values
of $r(A,B,B)$ attained by such pairs of sets.

We determine the minimum of the spectrum for $s\le2\cdot3^{n-1}$ and the maximum for $s\ge3^{n-1}$, and characterize when the upper bound $st$ is attained. 
We also show that for certain choices of $n,s$ and $t$, a block of consecutive values is always missing in the spectrum. 
Using subspaces and Sidon sets, we construct attainable values between the endpoints.
For Schur triples, we determine the minimum for
$t\le2\cdot3^{n-1}$ and the maximum for $t\ge3^{n-1}$.
We also extend our endpoint results to fields of odd prime power order.
\end{abstract}

\section{Introduction}\label{sec:intro}

We work in the vector space $V_n=\F_3^n$, where $n\ge 1$.
For subsets $A,B\subseteq\F_3^n$, an \emph{additive triple} is an ordered triple $(a,b,b')\in A\times B\times B$ satisfying $a+b=b'$.
Let
\[
 r(A,B,B):=\abs{\set{(a,b,b')\in A\times B\times B:a+b=b'}}.
\]
Equivalently, $r(A,B,B)$ is the number of pairs $(a,b)\in A\times B$ such that $a+b\in B$. 
The sets $A$ and $B$ need not be disjoint.
For each individual $a \in A$, we call $|B \cap (a+B)|$ the \emph{contribution} of $a$ to the value $r(A,B,B)$. It is clear that 
\begin{equation}\label{eq:contribution}
 |B \cap (a+B)| = |\{(b_1,b_2) \in B \times B: b_2-b_1 = a\}|.
\end{equation}
Therefore the contribution of $a$ can be determined by the number of times that $a$ occurs in the multiset of pairwise differences between elements of $B$.

When $A=B$, $r(A,B,B)$ is precisely the number of \emph{Schur triples} in $A$. And so the additive triple count generalises the Schur triple count by allowing the set $A$ to vary independently of $B$. In this paper, we study the possible values of $r(A,B,B)$ when the sizes of $A$ and $B$ are fixed. 
Let $N=3^n$.
For integers $0\le s,t\le N$, let
\[
 r_n(s,t):=\set{r(A,B,B):A,B\subseteq\F_3^n,\ \abs A=s,\ \abs B=t}.
\]
We call $r_n(s,t)$ the \emph{additive triple spectrum}.
Thus, if $r\in r_n(s,t)$, then there are sets $A$ and $B$ of sizes $s$ and $t$ with $r(A,B,B)=r$. 
If $r\notin r_n(s,t)$, then no such pair of sets exists.

Each triple can be determined by its first two entries, so $r(A,B,B)\le st$. 
It can also be determined by its last two entries, since $a=b'-b$, so $r(A,B,B)\le t^2$. 
Hence

\begin{equation}\label{eq:basic-range}
 r_n(s,t)\subseteq\set{0,1,\ldots,\min(st,t^2)}.
\end{equation}
But not every value inside the range can occur. 
See the following example.

\begin{example}\label{ex:two-three}
Let $n=2$, $s=2$, and $t=3$. It follows from \eqref{eq:basic-range} that $r_2(2,3) \subseteq \{0,1,\ldots,6\}$. 
In fact the spectrum is
\[
 r_2(2,3)=\set{0,1,2,3,4,6}.
\]
Let $a=(1,0)$ and $b=(0,1)$ be basis vectors of ${\mathbb F}_3^2$. We will demonstrate how each value of $r_2(2,3)$ can be obtained. Moreover, we will show that $5 \not\in r_2(2,3)$. Each pair of sets $A$ and $B$ in the following table attains the corresponding value of $r(A,B,B)$ in the left hand column.
\begin{center}
\begin{tabular}{ccl}
\toprule
$r(A,B,B)$ & $A$ & $B$ \\
\midrule
$0$ & $\set{b,2b}$ & $\set{0,a,2a}$ \\
$1$ & $\set{a,a+b}$ & $\set{0,a,b}$ \\
$2$ & $\set{a,2a}$ & $\set{0,a,b}$ \\
$3$ & $\set{0,b}$ & $\set{0,a,2a}$ \\
$4$ & $\set{0,a}$ & $\set{0,a,b}$ \\
$6$ & $\set{0,a}$ & $\set{0,a,2a}$ \\
\bottomrule
\end{tabular}
\end{center}
To see why $r(A,B,B) = 5$ is not attainable, first note that every affine line in $\F_3^2$ has the form $x+H$, where $H=\set{0,u,2u}$ and $u\ne0$.
When three points lie on the same line, they are \emph{collinear}.
Thus a $3$-subset of $\F_3^2$ is either an affine line or a set of non-collinear points.

Let $B=x+H$ be a line. Then
\[
 \abs{(B+y)\cap B}
 =
 \begin{cases}
 3,&y\in H,\\
 0,&y\notin H.
 \end{cases}
\]
It is immediate, that if $B = x + H$, then $r(A,B,B)$ takes the values $0$, $3$ and $6$ when $|A \cap H|$ is equal to $0$, $1$ and $2$, respectively. 

If $B$ is non-collinear, let $B=\set{x,x+u,x+v}$, where $u,v$ are linearly independent.

Let $\Delta(B)$ be the multi-set of pairwise differences 
\[
 \Delta(B):=\set{b-b':b,b'\in B}
 =\set{0,0,0,u,-u,v,-v,u-v,v-u}.
\]
It follows from the above and (\ref{eq:contribution}) that if $a = 0 \in A$ then the contribution of $a$ will be $3$, if $a \in A \setminus \{0\}$ and $a \in \Delta(B)$, then the contribution of $a$ to $r(A,B,B)$ will be $1$ and if $a \in A$ and $a \not\in \Delta(B)$ then the contribution of $a$ to $r(A,B,B)$ will be $0$.

Therefore, by summing all possibilities in both of the cases above, we find that $r_2(2,3) = \{0,1,2,3,4,6\}$.
Note that in both cases, there are no combinations of vectors that will produce an $r(A,B,B)$ value of $5$.

\end{example}

Similarly, it is not always possible to attain the upper and lower bounds of \eqref{eq:basic-range}.
This is demonstrated by the following example.

\begin{example}\label{ex:boundary}
In $\F_3^2$, we have
\[
 r_2(7,2)=\set{1,2,3,4},
 \qquad
 r_2(2,2)=\set{0,1,2,3}.
\]
Let $B=\{u,u+v\}$, where $v\ne0$. It follows that $\Delta (B) = \{0,0,v,2v\}$, therefore the elements $0,v,2v$ have contributions of $2,1,1$ to $r(A,B,B)$ respectively. 
As there are $9$ elements in $\mathbb{F}_3^2$, when $|A|=7$, $A$ contains at least one of these elements, so $r(A,B,B) = 0$ is impossible.

When $|A|=2$, then there are at most $2+1=3$ triples and only the zero shift has contribution $2$, so the value $t^2 = st = 4$ is not possible. 
\end{example}

For fixed $n,s,t$, which values occur in $r_n(s,t)$?
We denote 
\[
 f_n(s,t):=\min r_n(s,t),\qquad
 g_n(s,t):=\max r_n(s,t).
\]
Collectively, we call $f_n(s,t)$ and $g_n(s,t)$ the \textit{endpoints }of the spectrum $r_n(s,t)$.
If a value does not belong to the spectrum and lies strictly between the endpoints, we call it a \emph{hole}. 
Thus, $5$ is a hole in \Cref{ex:two-three}. 
The missing values in \Cref{ex:boundary} lie outside the endpoints, so they are not holes.

Both the endpoints and spectrum of possible values have been studied for Schur triples, which correspond to the case $A=B$.
A set containing no Schur triples is called \emph{sum-free}.
Diananda and Yap~\cite{DY} and Rhemtulla and Street~\cite{RS} studied the maximum size of a sum-free set for several classes of finite abelian groups.
Green and Ruzsa~\cite{GR} determined the maximum size for all finite abelian groups, and extended the results on counting the number of sum-free sets obtained by Lev, \L{}uczak and Schoen~\cite{LLS}.
Beyond the sum-free case, Huczynska, Mullen and Yucas~\cite{HMY} studied the possible numbers of Schur triples in a subset, focussing on groups of prime order.
Samotij and Sudakov~\cite{SS} determined the minimum number of Schur triples in several groups and size ranges. 
In particular, they obtain lower bounds on Schur triples for subsets of $\F_3^n$ slightly larger than the largest sum-free size.
Another direction has been to study the minimum number of additive tuples; see Chervak, Pikhurko and Staden~\cite{CPS} for groups of prime order.

When we allow $A$ and $B$ to vary, this leads to the additive triple spectrum. 
This study was initiated by Huczynska, Jedwab and Johnson~\cite{HJJ} in groups of odd prime order.
In their work, they determine the endpoints by applying Pollard's theorem~\cite{P}, and also show that every value between them is attained.
We also considered the generalization of Pollard's theorem to abelian groups ~\cite{GW}, but could not obtain improvements on the bounds proved here.
In a recent work, Huczynska, Rakhmonov and Yip~\cite[Theorem~1.5 and Corollary~4.2]{HRY} determine the endpoints of $r(A,B,B)$ in $\F_2^n$ for all given sizes of $A$ and $B$.
In contrast to the cyclic groups of prime order, in $\F_3^n$, subspaces and their cosets provide a different structure for constructions, and \Cref{ex:two-three} shows that holes can occur.

The paper will be organised as follows. In Section 2, we introduce notation and a series of preliminary results that will be used throughout this paper.

In Section 3, we begin our investigation into additive triples by determining the lower endpoint $f_n(s,t)$ for all sets in which $|A| = s \le 2N/3$.
We achieve this by viewing $r(A,B,B)$ as the sum of intersections of $B$ with its translates. 
A construction attaining this bound gives our first main result.

\begin{reptheorem}{thm:lower-main}
Let $N=3^n$. For $0\le s\le 2N/3$ and $0\le t\le N$,
\[
 f_n(s,t)=
 \begin{cases}
 0,&0\le t\le N/3,\\
 s(t-N/3),&N/3<t\le 2N/3,\\
 s(2t-N),&2N/3<t\le N.
 \end{cases}
\]
\end{reptheorem}

Next, we shift focus in Section 4 and look at $g_n(s,t)$. Specifically, we determine the values of $n,s$ and $t$ for which $g_n(s,t)=st$. The relation between the counts for $A$ and its complement also gives a formula for $g_n(s,t)$ when $s\ge N/3$. 

\begin{reptheorem}{thm:upper-main}
Let $N=3^n$. For $1\le s\le N$, let
$h_s=3^{\lceil\log_3s\rceil}$, the smallest power of $3$ that is
at least $s$. For $0\le t\le N$,
\[
 g_n(s,t)=st\quad\Longleftrightarrow\quad h_s\mid t.
\]
If $N/3\le s\le N$, then
\[
 g_n(s,t)=t^2-(N-s)\cdot \max\set{0,t-N/3,2t-N}.
\]
\end{reptheorem}

The cases $s\le3$ admit an exact formula.
Let $\varepsilon(t)=0$ when $3\mid t$ and
$\varepsilon(t)=1$ otherwise.
By \Cref{cor:small-s},
\[
 g_n(s,t)=st-(s-1)\varepsilon(t),
 \qquad 1\le s\le3.
\]
Our endpoint results also extend the range for the minimum number of Schur triples.
For $B\subseteq V_n$, let
\[
 \operatorname{ST}(B)
 =r(B,B,B)
 =\bigl|\{(x,y,z)\in B^3:x+y=z\}\bigr|.
\]
Note that triples are ordered i.e., $(x,y,z)$ and $(y,x,z)$ are counted separately when $x\ne y$, and their entries need not be distinct.
Let $N=3^n$ and $M=N/3$.
Samotij and Sudakov~\cite[Theorem~1.6]{SS} proved that
\[
 \operatorname{ST}(B)\ge Mu+u^2
 \qquad\text{when }|B|=M+u
\]
for $0\le u\le\delta M$, where $\delta>0$ is an absolute
constant, and provide sets attaining equality.
We extend the same formula throughout $0\le u\le N/3$.

\begin{reptheorem}{thm:schur-endpoints}
Let $n\ge1$, $N=3^n$, and $M=N/3$.
For every integer $0\le t\le2M$,

$$ \min_{\substack{B\subseteq V_n\\|B|=t}}
 \operatorname{ST}(B)
 =t\max\{0,t-M\}.$$

For every integer $M\le t\le N$,

$$ \max_{\substack{B\subseteq V_n\\|B|=t}}
 \operatorname{ST}(B)
 =
 \begin{cases}
 M^2+2(t-M)^2,&M\le t\le2M,\\[1mm]
 3t^2-3Nt+N^2,&2M\le t\le N.
 \end{cases}$$

\end{reptheorem}

In Section 5, we give bounds for $g_n(s,t)$ in the remaining range $4\le s<N/3$, $h_s\nmid t$. These also give bounds for $f_n(N-s,t)$. We prove the following theorem.

\begin{reptheorem}{thm:remaining-endpoint-bounds}
Let $N=3^n$, and let $s,t$ be integers with $4\le s<N/3$ and
$1\le t\le N$. Suppose that $h_s\nmid t$. Let $h_s=3^k$, $c_k(r)=\max\{0,r-h_s/3,2r-h_s\}$ and
$t=mh_s+r$, where $0<r<h_s$. Then
\[
 s(t-r)+r^2-(h_s-s)c_k(r)
 \le g_n(s,t)
 \le \min\{t^2,st-s+1\}.
\]
\end{reptheorem}

In Section 6, we extend these results to fields of odd prime power order. 

The endpoint results do not give us the complete picture. 
In Section 7, we proceed to show that the spectrum can contain long blocks of consecutive holes.
We prove the following result on \textit{holes}.

\begin{repproposition}{prop:top-blocks}
Let $N=3^n$, and let $s\ge3$ be a power of $3$.
Suppose that $s\mid t$ and $0<t<N$.
Then $g_n(s,t)=st$, and every value in
\[
 \{st-s+1,\ldots,st-1\}
\]
is a hole in $r_n(s,t)$.
\end{repproposition}

In Section 7, we also give constructions that attain values inside the spectrum, using subspaces and Sidon sets. 
We end with two open problems.

\section{Notation and preliminaries}\label{sec:prelim}

Let $V_n=\F_3^n$, where $n\ge1$,
and let $N=3^n$.
The parameters $s,t$ are integers with $0\le s,t\le N$.

For $a\in V_n$, translation by $a$ is the map $x\mapsto x+a$.
The corresponding translate of a set $B\subseteq V_n$ is
\[
 B+a:=\set{b+a:b\in B}.
\]
We will call $a$ the \emph{shift vector} or simply a \emph{shift}.

For a set $B\subseteq V_n$ of size $t$ and a vector $x\in V_n$,
define
\[
 N_x(B):=\abs{B\cap(B+x)}
 =\abs{\set{(b,b')\in B^2:b'-b=x}}.
\]
We will call $N_x(B)$ the \emph{overlap}. 
For $a\in A$, the overlap $N_a(B)$ is the number of additive triples with first entry $a$.
This gives us the following.
\begin{equation}\label{eq:r-sum-shifts-intro}
 r(A,B,B)=\sum_{a\in A}N_a(B).
\end{equation}
In particular,
\begin{equation}\label{eq:overlap-basics}
 N_0(B)=t,\qquad N_x(B)=N_{-x}(B),\qquad 0\le N_x(B)\le t.
\end{equation}
Translating $B$ does not change its overlaps:
\[
 N_x(B+y)=N_x(B)
 \qquad\text{for all }x,y\in V_n.
\]

We will also make use of the affine geometry of $V_n$.
We write $H\le V_n$ to show that $H$ is a linear subspace.
Its coset
\[
 y+H:=\set{y+h:h\in H}
\]
is an \emph{affine subspace}.
The dimension of $y+H$ is $\dim H$. It contains $3^{\dim H}$ points.
Two affine subspaces are called parallel if they are cosets
of the same linear subspace.
With these points and affine subspaces, $V_n$ is the affine
space $\operatorname{AG}(n,3)$.

For $a\ne0$, let
\[
 \angles a:=\set{\lambda a:\lambda\in\F_3}
 =\set{0,a,2a}.
\]
An \emph{affine line} is a coset
\[
 y+\angles a=\set{y,y+a,y+2a}.
\]
Its direction is the one-dimensional subspace $\angles a$.
Note that $a$ and $-a$ determine the same direction.
An \emph{affine hyperplane} is a coset of an
$(n-1)$-dimensional linear subspace. 
It has $N/3$ points.

We can also view $N_a(B)$ as the edge count in the subgraph induced by $B$. 
For $a\ne0$, let $\Gamma_a$ be the directed graph with vertex
set $V_n$ and directed edges
\[
 y\longrightarrow y+a
 \qquad (y\in V_n).
\]
This is the \emph{Cayley digraph}
$\operatorname{Cay}(V_n,\set{a})$ of the additive group $V_n$.
The subgraph induced by $B$ has $N_a(B)$ directed edges.

A directed $k$-cycle consists of distinct vertices
$v_0,\ldots,v_{k-1}$ and directed edges
$v_i\to v_{i+1}$, with indices $i$ modulo $k$.
For the zero shift, each vertex has a loop. The induced subgraph on $B$ therefore has $t$ loops.

We now prove our first identity.
\begin{lemma}\label{lem:total-overlap}
For every $B\subseteq V_n$ with $\abs B=t$,
\[
 \sum_{x\in V_n}N_x(B)=t^2.
\]
\end{lemma}
\begin{proof}
Each ordered pair $(b,b')\in B^2$ has a unique difference $x=b'-b$.
\end{proof}

Splitting the sum over $A$ and $V_n\setminus A$ gives a relation between the endpoints.

\begin{proposition}\label{prop:duality}
For $0\le s,t\le N$,
\[
 r_n(s,t)=\set{t^2-r:r\in r_n(N-s,t)}.
\]
Consequently,
\begin{equation}\label{eq:endpoint-duality}
 f_n(s,t)=t^2-g_n(N-s,t),\qquad
 g_n(s,t)=t^2-f_n(N-s,t).
\end{equation}
\end{proposition}
\begin{proof}
For $A,B\subseteq V_n$ with $|A|=s$ and $|B|=t$,
\Cref{lem:total-overlap} gives
\[
 \sum_{x\in A}N_x(B)
 +\sum_{x\in V_n\setminus A}N_x(B)
 =t^2.
\]
This proves the first statement of the proposition. 
Taking the minimum and maximum proves the second. 
\end{proof}

Thus, replacing $A$ by $V_n\setminus A$ replaces each value $r$ by $t^2-r$. 
This relates holes as well as endpoints: if $r$ is missing for $(s,t)$, then $t^2-r$ is missing for $(N-s,t)$.

\section{The lower endpoint}\label{sec:lower}

We first minimize the overlap $N_a(B)$ for one fixed nonzero
shift vector $a$.
We then construct a set $B$ for which $2N/3$ different
shifts attain the same minimum.
For $s\le2N/3$, choosing $A$ from these shifts attains the lower bound.

Fix $a\in V_n\setminus\set{0}$.
Let $\Gamma_a$ be the directed graph defined as before.
The cosets of $\angles a$ partition $V_n$ into $M=N/3$
affine lines.
On each line $L=\set{y,y+a,y+2a}$, translation by $a$
gives the directed $3$-cycle
\[
 y\longrightarrow y+a\longrightarrow y+2a\longrightarrow y.
\]
Thus the affine lines parallel to $\angles a$ are precisely
the vertex sets of the directed cycles of $\Gamma_a$.

Let $k_L=\abs{B\cap L}$.
The number of directed edges on $L$ with both endpoints in $B$
depends only on $k_L$.
Denote this number by
\[
 \phi(k_L)
 :=\abs{\set{b\in B\cap L:b+a\in B\cap L}}.
\]
If there is zero or one vertex in $B \cap L$, that gives no edges. 
Two vertices give one edge, and all three vertices give three edges.
Hence,
\[
 \phi(0)=\phi(1)=0,\qquad \phi(2)=1,\qquad \phi(3)=3.
\]
And when we sum over the cosets of $\angles a$, this gives us
\begin{equation}\label{eq:cycle-rule}
 N_a(B)=\sum_L\phi(k_L),\qquad \sum_Lk_L=t.
\end{equation}

The following example illustrates this.

\begin{example}\label{ex:translation-lines}
Let $V_2=\F_3^2$ and $a=(1,0)$.
The cosets of $\angles a$ are the three parallel lines 
\[
 L_j=\set{(0,j),(1,j),(2,j)},
 \qquad j\in\F_3.
\]
Let
\[
 B=\set{(0,0),(1,0),(0,1),(0,2)}.
\]
Then
\[
 |B\cap L_0|=2,\qquad
 |B\cap L_1|=1,\qquad
 |B\cap L_2|=1.
\]
Only $L_0$ contains a directed edge with both endpoints in $B$,
namely $(0,0)\to(1,0)$.
Thus, when $A=\set{a}$,
\[
 r(A,B,B)=N_a(B)
 =\phi(2)+\phi(1)+\phi(1)=1.
\]

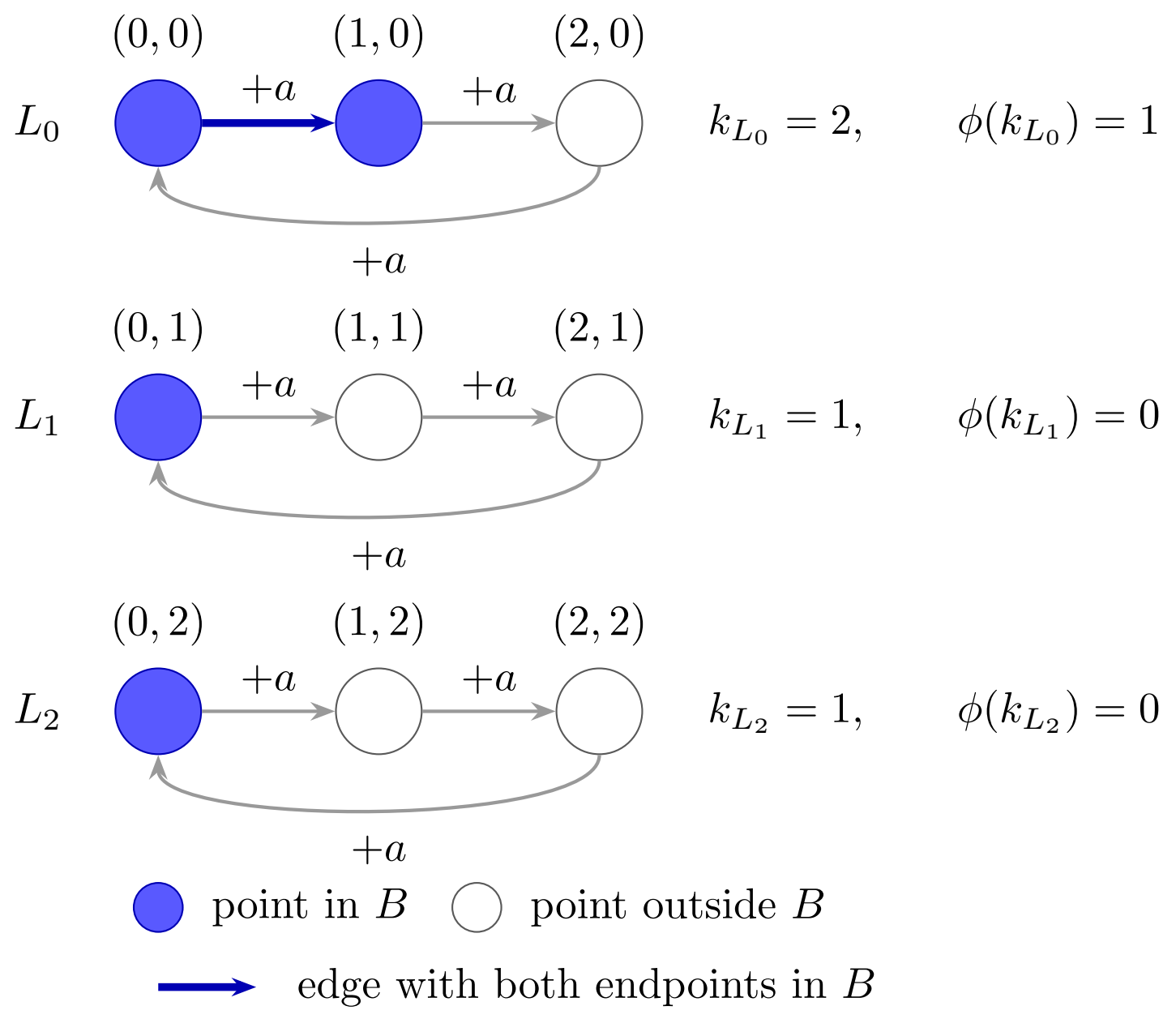
\begin{figure}[htbp]
\centering
\begin{tikzpicture}[
  vertex/.style={
    circle, draw=black!65, fill=white,
    minimum size=7mm, inner sep=0pt
  },
  selected/.style={
    vertex, draw=blue!70!black, fill=blue!65
  },
  translationarrow/.style={
    -{Stealth[length=2mm]},
    draw=black!40, line width=0.8pt
  },
  counted/.style={
    -{Stealth[length=2mm]},
    draw=blue!70!black, line width=1.8pt
  },
  every node/.style={font=\small}
]

\foreach \j/\k/\f in {0/2/1,1/1/0,2/1/0}{
  \node[anchor=east] at (-0.65,{-2.4*\j})
    {$L_{\j}$};

  \foreach \i in {0,1,2}{
    \ifnum\i<\k\relax
      \node[selected] (v-\j-\i)
        at ({1.8*\i},{-2.4*\j}) {};
    \else
      \node[vertex] (v-\j-\i)
        at ({1.8*\i},{-2.4*\j}) {};
    \fi

    \node[above=4mm] at (v-\j-\i)
      {$(\i,\j)$};
  }

  \ifnum\j=0\relax
    \draw[counted]
      (v-\j-0.east) --
      node[midway,above] {$+a$}
      (v-\j-1.west);
  \else
    \draw[translationarrow]
      (v-\j-0.east) --
      node[midway,above] {$+a$}
      (v-\j-1.west);
  \fi

  \draw[translationarrow]
    (v-\j-1.east) --
    node[midway,above] {$+a$}
    (v-\j-2.west);

  \draw[translationarrow]
    (v-\j-2.south)
    .. controls (3.6,{-2.4*\j-0.95})
            and (0,{-2.4*\j-0.95}) ..
    node[midway,below=2pt] {$+a$}
    (v-\j-0.south);

  \node[anchor=west] at (4.35,{-2.4*\j})
    {$k_{L_{\j}}=\k,\qquad \phi(k_{L_{\j}})=\f$};
}

\node[selected,minimum size=4mm] at (0,-6.4) {};
\node[anchor=west] at (0.3,-6.4)
  {point in $B$};

\node[vertex,minimum size=4mm] at (2.6,-6.4) {};
\node[anchor=west] at (2.9,-6.4)
  {point outside $B$};

\draw[counted] (0,-7.05) -- (0.8,-7.05);
\node[anchor=west] at (1,-7.05)
  {edge with both endpoints in $B$};

\end{tikzpicture}
\caption{Translation by $a=(1,0)$ on the three parallel lines
of $\operatorname{AG}(2,3)$.
Blue vertices belong to $B$.
The blue directed edge is the only edge with both endpoints in $B$.}
\label{fig:three-line-example}
\end{figure}
\end{example}

Next, we minimize $N_a(B)$ over all sets of size $t$.
Define
\[
 c_n(t):=
 \begin{cases}
  0, & 0\le t\le N/3,\\
  t-N/3, & N/3<t\le 2N/3,\\
  2t-N, & 2N/3<t\le N.
 \end{cases}
\]
\begin{lemma}\label{lem:one-shift-min}
For every nonzero $a\in V_n$ and $0\le t\le N$,
\[
 \min_{\abs B=t}N_a(B)=c_n(t).
\]

\end{lemma}

\begin{proof}
For $k\in\{0,1,2,3\}$, we have
\[
 \phi(k)\ge0,\qquad
 \phi(k)\ge k-1,\qquad
 \phi(k)\ge2k-3.
\]
Since there are $N/3$ lines and $\sum_L k_L=t$,
summing these inequalities gives
\[
 N_a(B)\ge\max\{0,t-N/3,2t-N\}=c_n(t).
\]

If $t\le N/3$, choose one point on each of $t$ distinct
lines. This gives $N_a(B)=0$.
If $N/3<t\le2N/3$, place one
point on $2N/3-t$ lines and two points on $t-N/3$ lines. This gives $N_a(B)=t-N/3$.
If $2N/3<t\le N$, place two
points on $N-t$ lines and three points on $t-2N/3$ lines. This gives $N_a(B)=2t-N$.
\end{proof}

Since $N_0(B)=t\ge c_n(t)$, every shift has overlap at least $c_n(t)$.
Therefore, for all $A,B$ of sizes $s,t$,
\begin{equation}\label{eq:universal-lower}
 r(A,B,B)\ge s\,c_n(t).
\end{equation}

The bound in \Cref{lem:one-shift-min} holds for every set $B$.
To attain the lower bound for $r(A,B,B)$, we need one set $B$
for which all shifts chosen in $A$ attain this minimum.
We now construct a set with this property for $2N/3$
different shift vectors at once.
The construction is inspired by Samotij and Sudakov~\cite[Theorem~1.6]{SS},
who take a full affine hyperplane and part of a parallel one.

Let
\[
 V_n=U\oplus\F_3,\qquad
 U=\F_3^{n-1},\qquad |U|=M=N/3.
\]
When $n=1$, $U$ is the zero-dimensional vector space,
which has one element.
For $0\le t<N$, let $t=qM+r$ with $0\le q\le2$ and $0\le r<M$.
Choose any $T\subseteq U$ with $\abs T=r$.
Define
\begin{equation}\label{eq:layered-B}
 B^\downarrow(t):=
 \begin{cases}
 T\times\set{0},&q=0,\\
 (U\times\set{0})\cup(T\times\set{1}),&q=1,\\
 (U\times\set{0,1})\cup(T\times\set{2}),&q=2.
 \end{cases}
\end{equation}
If $t=N$, then let $B^\downarrow(N)=V_n$.
In each case the set has size $t$.
For $i\in\F_3$, call $U\times\set{i}$ the $i$th \emph{layer}.
The three layers are parallel affine hyperplanes.
We fill the layers in order. At most one layer is partly filled.

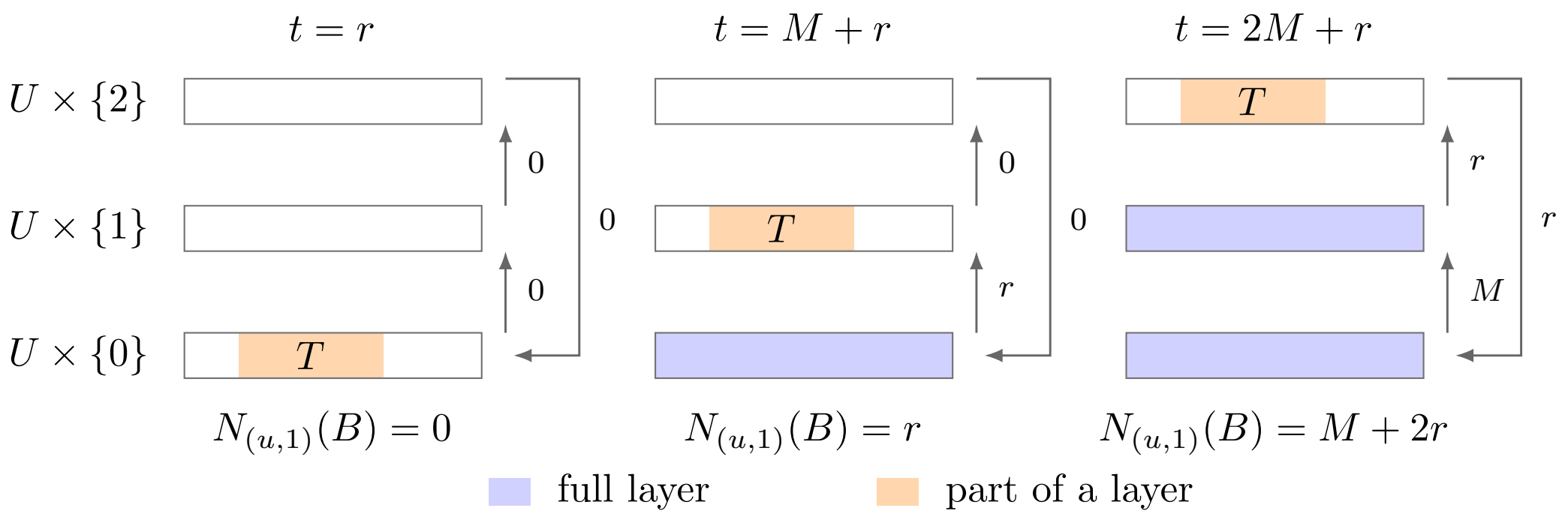
\begin{figure}[htbp]
\centering
\begin{tikzpicture}[x=1cm,y=1cm,
  every node/.style={font=\small},
  shiftarrow/.style={-{Latex[length=1.6mm]},draw=black!60,line width=0.55pt}]
\foreach \i in {0,1,2}{
  \node[anchor=east] at (-0.18,{1.12*\i+0.20}) {$U\times\set{\i}$};
}
\foreach \q in {0,1,2}{
 \begin{scope}[xshift={4.15*\q cm}]
  \foreach \i in {0,1,2}{
   \ifnum\i<\q
    \fill[blue!18] (0,{1.12*\i}) rectangle (2.62,{1.12*\i+0.40});
   \else
    \ifnum\i=\q
     \fill[orange!32] (0.48,{1.12*\i}) rectangle (1.75,{1.12*\i+0.40});
     \node at (1.11,{1.12*\i+0.20}) {$T$};
    \fi
   \fi
   \draw[black!55] (0,{1.12*\i}) rectangle (2.62,{1.12*\i+0.40});
  }
  \draw[shiftarrow] (2.83,0.40) -- (2.83,1.12);
  \draw[shiftarrow] (2.83,1.52) -- (2.83,2.24);
  \draw[shiftarrow] (2.83,2.64) -- (3.48,2.64) -- (3.48,0.20) -- (2.90,0.20);
 \end{scope}
}
\node at (1.30,3.08) {$t=r$};
\node at (5.45,3.08) {$t=M+r$};
\node at (9.60,3.08) {$t=2M+r$};
\node[anchor=west,font=\scriptsize] at (2.89,0.78) {$0$};
\node[anchor=west,font=\scriptsize] at (2.89,1.90) {$0$};
\node[anchor=west,font=\scriptsize] at (3.52,1.40) {$0$};
\node[anchor=west,font=\scriptsize] at (7.04,0.78) {$r$};
\node[anchor=west,font=\scriptsize] at (7.04,1.90) {$0$};
\node[anchor=west,font=\scriptsize] at (7.67,1.40) {$0$};
\node[anchor=west,font=\scriptsize] at (11.19,0.78) {$M$};
\node[anchor=west,font=\scriptsize] at (11.19,1.90) {$r$};
\node[anchor=west,font=\scriptsize] at (11.82,1.40) {$r$};
\node at (1.30,-0.48) {$N_{(u,1)}(B)=0$};
\node at (5.45,-0.48) {$N_{(u,1)}(B)=r$};
\node at (9.60,-0.48) {$N_{(u,1)}(B)=M+2r$};
\fill[blue!18] (2.68,-1.12) rectangle (3.05,-0.88);
\node[anchor=west] at (3.15,-1.00) {full layer};
\fill[orange!32] (6.10,-1.12) rectangle (6.47,-0.88);
\node[anchor=west] at (6.57,-1.00) {part of a layer};
\end{tikzpicture}
\caption{
The three layers of $B^\downarrow(t)$, with $0\le r<M$.
Blue layers are full, orange regions represent an arbitrary $r$-point set $T$, and unshaded regions are empty. 
A shift $(u,1)$ sends layer $0$ to layer $1$, layer $1$ to layer $2$, and layer $2$ to layer $0$. 
The arrow labels give the contributions of the three overlaps.
They do not depend on $u$. 
At $r=0$, the orange region is empty. At $t=N$, all three layers are full.
}
\label{fig:layers}
\end{figure}

A shift vector $(u,1)$, where $u\in U$, sends
$(v,i)$ to $(v+u,i+1)$, with the second coordinate is in $\F_3$.
And so $(u,1)$ sends layer $0$ to layer $1$, layer $1$ to layer $2$, and layer $2$ to layer $0$.
\Cref{fig:layers} shows the resulting contributions from the overlaps.

The shifts used in the construction are precisely the vectors
outside $U\times\set{0}$.
Let
\[
 W:=U\times\set{1,2}.
\]
Then $|W|=2N/3$ and $0\notin W$.

\begin{lemma}\label{lem:many-shifts}
For $B=B^\downarrow(t)$ and every $a\in W$,
\[
 N_a(B)=c_n(t).
\]
\end{lemma}
\begin{proof}
For $i\in\F_3$, let $U_i:=\set{v\in U:(v,i)\in B}$.
Then
\[
 B=\bigcup_{i=0}^2(U_i\times\set{i}).
\]
For $a=(u,1)$,
\[
 N_a(B)=\abs{U_1\cap(U_0+u)}
       +\abs{U_2\cap(U_1+u)}
       +\abs{U_0\cap(U_2+u)}.
\]
The layers in \eqref{eq:layered-B} give $0$, $r$, and
$M+2r$. These are $c_n(t)$ in the corresponding ranges.
At $t=N$, the overlap is $N=c_n(N)$.
For $a=(u,2)$, use $N_a(B)=N_{-a}(B)$, since $-a=(-u,1)$.
\end{proof}

Choosing $A$ (shifts) in such a way that minimizes the number of triples gives us our first main result of the paper. 

\begin{theorem}\label{thm:lower-main}
For $0\le s\le2N/3$ and $0\le t\le N$,
\[
 f_n(s,t)=s\,c_n(t).
\]
\end{theorem}
\begin{proof}
The lower bound is \eqref{eq:universal-lower}.
For equality, take $B=B^\downarrow(t)$ and choose any $s$-point subset
$A\subseteq W$. Each shift in $A$ has overlap $c_n(t)$ by
\Cref{lem:many-shifts}, so $r(A,B,B)=s\,c_n(t)$.
\end{proof}

In particular, the lower endpoint is zero for $t\le N/3$ and
$s\le2N/3$. 
For $t>N/3$, every nonzero shift has positive overlap, so zero is excluded from the spectrum whenever $s>0$.

\section{The upper endpoint}\label{sec:upper}

Recall from our earlier discussion that $r_n(s,t)\subseteq\set{0,1,\ldots,\min(st,t^2)}$. 
So we first ask when the bound $st$ is attained.
When does a shift maximize the overlap?
For this, we first define the following parameter, which we will use throughout the section.

For $m\ge1$, let
\begin{equation}\label{eq:hs}
 h_m:=3^{\lceil\log_3m\rceil}.
\end{equation}
$h_m$ is the smallest power of $3$, that is at least $m$.
We use $h_m$ only for $m\le N$. A subspace of this size exists in $V_n$.
The cases $s=0$ or $t=0$ give $r_n(s,t)=\set{0}$.
If $s=N$, then $r_n(N,t)=\set{t^2}$; if $t=N$, then
$r_n(s,N)=\set{sN}$.

A shift $x$ has overlap $t$ exactly when $B+x=B$.
For $B\subseteq V_n$, define its \emph{stabilizer} by
\[
 \Stab(B):=\set{x\in V_n:B+x=B}.
\]
We say that $x\in V_n$ \emph{stabilizes} $B$ if $B+x=B$.
$\Stab(B)$ is an additive subgroup of $V_n$, and hence a subspace over $\F_3$.
The set $B$ is a union of cosets of $\Stab(B)$.
These statements also hold for $B=\varnothing$.

We begin by proving the relationship between the stabilizer and maximum overlap.
\begin{lemma}\label{lem:stabilizer}
For every $x\in V_n$,
\[
 N_x(B)=t\quad\Longleftrightarrow\quad x\in\Stab(B).
\]
Consequently,
\[
 r(A,B,B)=st\quad\Longleftrightarrow\quad A\subseteq\Stab(B).
\]
\end{lemma}
\begin{proof}
If $N_x(B)=t$, then $|B\cap(B+x)|=|B|$, and so $B\subseteq B+x$.
Since $|B+x|=|B|$, this gives $B+x=B$.
Conversely, if $B+x=B$, then $N_x(B)=|B|=t$.
Next, by \eqref{eq:r-sum-shifts-intro}, $r(A,B,B)=\sum_{a\in A}N_a(B)$, where $N_a(B) \le t, \; \forall a \in A$.
Hence $r(A,B,B)=st$ exactly when $N_a(B)=t$ for every $a\in A$.
By the first part, this is equivalent to $A\subseteq\Stab(B)$.
\end{proof}

Now we prove our first result on $g_n(s,t)$ for certain regimes of $s$ and $t$.
\begin{theorem}\label{thm:upper-main}
For $1\le s\le N$ and $0\le t\le N$,
\[
 g_n(s,t)=st\quad\Longleftrightarrow\quad h_s\mid t.
\]
For $N/3\le s\le N$ and $0\le t\le N$,
\[
 g_n(s,t)=t^2-(N-s)c_n(t).
\]
\end{theorem}
\begin{proof}
Suppose $r(A,B,B)=st$, and let $H=\Stab(B)$.
By \Cref{lem:stabilizer}, $A\subseteq H$, so $\abs H\ge s$.
Since $\abs H$ is a power of $3$, we have $h_s\mid\abs H$.
Since $B$ is a union of $H$-cosets, $\abs H\mid t$.
Therefore $h_s\mid t$.

Conversely, if $h_s\mid t$, choose a subspace $H$ of size $h_s$.
Let $B$ be a union of $t/h_s$ cosets of $H$, and choose an $s$-point
subset $A\subseteq H$.
Every element of $A$ stabilizes $B$, so $r(A,B,B)=st$.
This also covers $t=0$.

For the second assertion, $s\ge N/3$ implies $N-s\le2N/3$.
Apply \Cref{thm:lower-main} and \eqref{eq:endpoint-duality}.
\end{proof}

This gives us the following formula for the lower-endpoint.

\begin{corollary}\label{cor:large-lower}
Let $u=N-s$. If $u=0$, then $f_n(s,t)=t^2$.
If $u\ge1$ and $h_u\mid t$, then
\[
 f_n(s,t)=t^2-ut.
\]
\end{corollary}
\begin{proof}
From \Cref{thm:upper-main}, we have $h_u \mid t \iff g_n(u,t)=ut$. 
Now apply \eqref{eq:endpoint-duality}. 
\end{proof}
Next we provide (sufficient) conditions for attaining the universal upper bound of $t^2$.

\begin{lemma}\label{lem:t-squared}
If $1\le t\le N$ and $h_t\le s\le N$, then
\[
 g_n(s,t)=t^2.
\]
\end{lemma}

\begin{proof}
Choose a subspace $L\le V_n$ of size $h_t$,
a set $B\subseteq L$ of size $t$, and a set
$A\subseteq V_n$ of size $s$ containing $L$.
And so we have $B \subseteq L \subseteq A$.
Every ordered pair $(b,b')\in B^2$ has its difference
$b'-b$ in $L\subseteq A$, so $r(A,B,B)=t^2$.
Since $r(A,B,B)\le t^2$ for every pair $A,B$, the maximum is $t^2$.
\end{proof}
Note: To count all $t^2$ pairs, every difference $b-b'$ must belong to the set $A$. 
And hence $B-B\subseteq A$. Combining this with ~\cite[Theorem~7]{L}, the condition in \Cref{lem:t-squared} is also necessary.

For a nonzero shift, the largest possible overlap is the following.

\begin{lemma}\label{lem:one-shift-max}
For $a\ne0$ and $0\le t\le N$,
\[
 \max_{\abs B=t}N_a(B)=
 \begin{cases}
 t,&3\mid t,\\
 t-1,&3\nmid t.
 \end{cases}
\]
\end{lemma}
\begin{proof}
By \Cref{lem:stabilizer}, $N_a(B)=t \iff a\in\Stab(B)$. 
In this case $\angles a\subseteq\Stab(B)$, so
$3\mid\abs{\Stab(B)}$. Since $B$ is a union of cosets of $\Stab(B)$,
we have $\abs{\Stab(B)}\mid t$, and hence $3\mid t$.
Thus $N_a(B)\le t-1$ when $3\nmid t$.

To attain the upper bound, construct $B$ as follows.
Take the union of $\lfloor t/3\rfloor$ cosets of $\angles a$.
If $3\mid t$, this gives $N_a(B)=t$.
Otherwise, add the remaining one or two points from one further coset.
That coset contributes zero to the overlap if it contains one point of $B$
and one if it contains two. Thus in either case $N_a(B)=t-1$.
\end{proof}

And now we can determine $g_n(s,t)$ for $1\le s\le3$.

\begin{corollary}\label{cor:small-s}
Let $\varepsilon(t)=0$ if $3\mid t$ and $\varepsilon(t)=1$ otherwise.
For $1\le s\le3$,
\[
 g_n(s,t)=st-(s-1)\varepsilon(t).
\]

\end{corollary}
\begin{proof}
We have $h_1=1$ and $h_2=h_3=3$.
Now if $3\mid t$, then $h_s\mid t$. 
Thus \Cref{thm:upper-main} gives $g_n(s,t)=st$.

Suppose $3\nmid t$. By \Cref{lem:one-shift-max}, every nonzero shift
has overlap at most $t-1$. The zero shift has overlap $t$ and can
occur at most once in the sum over $A$. Hence
\[
 r(A,B,B)\le t+(s-1)(t-1)=st-(s-1).
\]

For equality, fix $a\ne0$ and choose $B$ of size $t$ with
$N_a(B)=t-1$, as in \Cref{lem:one-shift-max}.
Since $N_{-a}(B)=N_a(B)$, take
\[
 A=
 \begin{cases}
 \{0\},&s=1,\\
 \{0,a\},&s=2,\\
 \{0,a,-a\},&s=3.
 \end{cases}
\]
In each case, the zero shift contributes $t$ and the other shifts contribute $t-1$ each.
\end{proof}
Together with \eqref{eq:endpoint-duality} and \Cref{thm:lower-main}, this determines
both endpoints when $s\le3$ or $N-s\le3$.

\subsection{Schur triples}
\label{subsec:schur-triples}

Our theorems on endpoints can be applied in the setting where $A=B$.

\begin{theorem}
\label{thm:schur-endpoints}
Let $n\ge1$, $N=3^n$, and $M=N/3$.
For every integer $0\le t\le2M$,
\begin{equation}\label{eq:schur-minimum}
 \min_{\substack{B\subseteq V_n\\|B|=t}}
 \operatorname{ST}(B)
 =t\max\{0,t-M\}.
\end{equation}
For every integer $M\le t\le N$,
\begin{equation}\label{eq:schur-maximum}
 \max_{\substack{B\subseteq V_n\\|B|=t}}
 \operatorname{ST}(B)
 =
 \begin{cases}
 M^2+2(t-M)^2,&M\le t\le2M,\\[1mm]
 3t^2-3Nt+N^2,&2M\le t\le N.
 \end{cases}
\end{equation}
\end{theorem}

\begin{proof}
Let $B\subseteq V_n$ have size $t$.
Recall that
\[
 c_n(t)=\max\{0,t-M,2t-N\}.
\]
By \Cref{lem:one-shift-min} and the identity $N_0(B)=t$,
every shift satisfies $N_a(B)\ge c_n(t)$.
Since
\[
 \operatorname{ST}(B)=\sum_{a\in B}N_a(B),
 \qquad
 \sum_{a\in V_n}N_a(B)=t^2,
\]
we obtain
\begin{equation}\label{eq:schur-universal-bounds}
 t\,c_n(t)
 \le \operatorname{ST}(B)
 \le t^2-(N-t)c_n(t).
\end{equation}

Now identify $V_n=H\times\F_3$, where $H=\F_3^{n-1}$ and
$|H|=M$.

For $0\le t\le M$, choose a $t$-point subset $B$ of $H\times\{1\}$.
The sum of any two points of $B$ lies in $H\times\{2\}$, so
$\operatorname{ST}(B)=0$. Thus the minimum is zero.

Suppose $M\le t\le2M$, and let $u=t-M$.
Then $c_n(t)=u$, and so 
\[
 \operatorname{ST}(B) \ge t\,c_n(t)=(M+u)u=Mu+u^2.
\]
Now we construct $B$ which attains this lower bound.
Choose any $T\subseteq H$ with $|T|=u$, and let
\[
 B=(H\times\{1\})\cup(T\times\{2\}).
\]
\begin{enumerate}
    \item For each $z\in T\times\{2\}$ and each $x\in H\times\{1\}$, the point $y=z-x$ lies in $H\times\{1\}$. This gives $Mu$ triples.
    \item If $x,y\in T\times\{2\}$, then $x+y\in H\times\{1\}$, giving another $u^2$ triples. 
    \item If one of $x,y$ lies in $H\times\{1\}$ and the other lies in $T\times\{2\}$, then $x+y\in H\times\{0\}$, which is outside $B$.
\end{enumerate}
Thus $\operatorname{ST}(B)=Mu+u^2$.
This attains the lower bound and proves \eqref{eq:schur-minimum}.

For the maximum in the same range, 
\[
 \begin{aligned}
 \operatorname{ST}(B) \le t^2-(N-t)c_n(t)
 &=(M+u)^2-(2M-u)u\\
 &=M^2+2u^2.
 \end{aligned}
\]
To attain the upper bound, choose $T\subseteq H$ with $|T|=u$ and let
\[
 B=(H\times\{0\})\cup(T\times\{1\}).
\]
\begin{enumerate}
 \item If $x,y\in H\times\{0\}$, then $x+y\in H\times\{0\}$.
 This gives $M^2$ triples.
 \item For each $y,z\in T\times\{1\}$, the point $x=z-y$ lies in
 $H\times\{0\}$. This gives $u^2$ triples $(x,y,z)$ and another
 $u^2$ triples $(y,x,z)$.
 \item If $x,y\in T\times\{1\}$, then $x+y\in H\times\{2\}$,
 outside $B$.
\end{enumerate}
Thus $\operatorname{ST}(B)=M^2+2u^2$, which proves the first case
of \eqref{eq:schur-maximum}.

Finally, suppose $2M\le t\le N$. Then $c_n(t)=2t-N$.
Choose $C\subseteq H\times\{1\}$ with $|C|=N-t$ and let
$B=V_n\setminus C$. For every $a\in C$, we have
$C+a\subseteq H\times\{2\}$, so $C$ and $C+a$ are disjoint.
Since $B\cap(B+a)
=(V_n\setminus C)\cap\bigl(V_n\setminus(C+a)\bigr)
=V_n\setminus\bigl(C\cup(C+a)\bigr)$,
then we have
\[
 N_a(B)=N-2|C|=2t-N.
\]
Subtracting the contributions from $a\in C$ gives us
\[
 \operatorname{ST}(B)
 =t^2-\sum_{a\in C}N_a(B)
 =t^2-(N-t)(2t-N)
 =3t^2-3Nt+N^2.
\]
This attains the upper bound in \eqref{eq:schur-universal-bounds}
and proves the second case of \eqref{eq:schur-maximum}.
\end{proof}

\section{Remaining range}\label{sec:bounds}
In \Cref{sec:lower}, we determined $f_n(s,t)$ for $s\le 2N/3$
(\Cref{thm:lower-main}). In \Cref{sec:upper}, we determined $g_n(s,t)$
for $s\ge N/3$ and showed that $g_n(s,t)=st$ exactly when $h_s\mid t$
(\Cref{thm:upper-main}). We also determined $g_n(s,t)$ for
$1\le s\le3$ in \Cref{cor:small-s}.
We now consider the range
\begin{equation}\label{eq:remaining}
 4\le s<N/3,\qquad h_s\nmid t.
\end{equation}
By \eqref{eq:endpoint-duality}, we know that the bounds on $g_n(s,t)$ also give bounds for $f_n(N-s,t)$.
We first give a lower bound and then an upper bound on $g_n(s,t)$.

Let $h=3^k$ and let
\[
 c_k(r)=\max\{0,r-h/3,2r-h\}.
\]
This is the one-shift minimum for an $r$-point set inside a $k$-dimensional space.
The following proposition gives a lower bound for $g_n(s,t)$.
We do this by constructing one pair $A,B$ of sizes $s,t$. 
The idea is to choose $A$ inside a subspace $H$, so that every shift in $A$ acts separately on each coset of $H$.

\begin{proposition}\label{prop:flat-bound}
Let $4\le s<N/3$ and $1\le t\le N$, and suppose that $h_s\nmid t$.
Let $h_s=3^k$ and $t=mh_s+r$, where $0<r<h_s$. Then
\begin{equation}\label{eq:flat-bound}
 g_n(s,t)\ge s(t-r)+r^2-(h_s-s)c_k(r).
\end{equation}
\end{proposition}

\begin{proof}
Choose a subspace $H\le V_n$ of size $h_s$. Construct $B$ from
$m$ full cosets of $H$ and a set $x+T$ in one further coset $x+H$,
where $T\subseteq H$ and $|T|=r$. Such a coset exists because
$mh_s<t\le N$. Take $A\subseteq H$ with $|A|=s$.

Every shift in $A$ preserves each coset of $H$. Thus each full coset
contributes $sh_s$ triples, and together the full cosets contribute
$smh_s=s(t-r)$.

For $a\in A$ and $b\in T$, shifting $x+b$ by $a$ gives $x+(a+b)$.
This point lies in $x+T$ exactly when $a+b\in T$.
So the partial coset contributes exactly $r(A,T,T)$.

Since $h_s$ is the smallest power of $3$ at least $s$, we have
$s>h_s/3$. Applying \Cref{thm:upper-main} inside $H$, we can choose
$A$ and $T$ so that $r(A,T,T)=r^2-(h_s-s)c_k(r)$.
\end{proof}

The following proposition gives the upper bound on $g_n(s,t)$.

\begin{proposition}\label{prop:stabilizer-bound}
Let $4\le s<N/3$ and $1\le t\le N$, and suppose that
$h_s\nmid t$. Then
\begin{equation}\label{eq:simple-upper}
 g_n(s,t)\le\min\{t^2,st-s+1\}.
\end{equation}
\end{proposition}

\begin{proof}
For a pair $A,B$ of sizes $s,t$, let $H=\Stab(B)$
and $|H|=\ell$.
Since $H$ is a subspace and $B$ is a union of $H$-cosets,
$\ell$ is a power of $3$ dividing $t$.
If $\ell\ge h_s$, then $h_s\mid\ell$, because both
are powers of $3$.
Since $\ell\mid t$, we would have $h_s\mid t$,
contradicting the assumption.
So $\ell<h_s$, and $\ell\le h_s/3<s$.

Both $B$ and $B+a$ are unions of $H$-cosets, so their
intersection is also a union of $H$-cosets.
Hence $N_a(B)$ is divisible by $\ell$.
For $a\notin H$, we have $N_a(B)<t$, and therefore
$N_a(B)\le t-\ell$.

Let $m=|A\cap H|$.
The $m$ shifts in $H$ contribute $t$ each, and the other
$s-m$ shifts contribute at most $t-\ell$ each.
Since $m\le\ell$,
\[
 \begin{aligned}
 r(A,B,B)
 &\le mt+(s-m)(t-\ell)\\
 &=st-(s-m)\ell\\
 &\le st-(s-\ell)\ell\\
 &\le st-s+1,
 \end{aligned}
\]
where the last inequality follows from $(s-\ell)\ell-(s-1)  =(\ell-1)(s-\ell-1)\ge0.$
Combining this with the universal bound $r(A,B,B)\le t^2$
proves the result.
\end{proof}

We now collect and present the bounds from the two propositions.

\begin{theorem}\label{thm:remaining-endpoint-bounds}
Let $N=3^n$, and let $s,t$ be integers with $4\le s<N/3$ and
$1\le t\le N$. Suppose that $h_s\nmid t$. Let $h_s=3^k$ and
$t=mh_s+r$, where $0<r<h_s$. Then
\[
 s(t-r)+r^2-(h_s-s)c_k(r)
 \le g_n(s,t)
 \le \min\{t^2,st-s+1\}.
\]
\end{theorem}
\begin{proof}
The lower bound follows from \Cref{prop:flat-bound}, and the upper
bound follows from \Cref{prop:stabilizer-bound}.
\end{proof}

Using \eqref{eq:endpoint-duality}, we get bounds for $f_n(N-s,t)$.

\begin{remark}\label{rem:remaining-bounds-sharp}
Let $h_s=3^k$ and $t=mh_s+r$ as in
\Cref{thm:remaining-endpoint-bounds}.

If $r\in\{1,h_s-1\}$, then the lower bound is
$st-s+1$, since $c_k(1)=0$ and
$c_k(h_s-1)=h_s-2$.
Thus $g_n(s,t)=st-s+1$ when
$t\equiv\pm1\pmod{h_s}$.
\end{remark}

The bounds also coincide in the following family.

\begin{remark}\label{rem:distinct-bounds-exact}
Let $n\ge3$, $s=4$, and $t=9m+6$, where
$0\le m<3^{n-2}$.
Here $h_s=9$ and $c_2(6)=3$, so
\Cref{thm:remaining-endpoint-bounds} gives $4t-3\le g_n(4,t)\le4t-3$.
Thus $g_n(4,t)=4t-3$.
\end{remark}

\section{Extending to fields of odd prime-power order}\label{sec:odd-prime}

In this section, we will extend some of our results to $\F_p^n$, where $p$ is an odd prime. 
Our theorem also applies to $\mathbb F_{p^e}^{\,n}$ by viewing $\mathbb F_{p^e}^{\,n}\cong\mathbb F_p^{\,en}$ as vector spaces over $\mathbb F_p$. 
This preserves addition, so it preserves all the additive triple counts. 

Let $p \ge 3$ be prime, and let $V_{p,n}=\F_p^n, \; N_p=p^n, \; M=p^{n-1},\; \text{and } p=2q+1$.
Define $r_{p,n}(s,t)$, $f_{p,n}(s,t)$, and $g_{p,n}(s,t)$ as in the previous sections.
Let $c_{p,n}(t):=\max\set{0,t-qM,2t-N_p}$ and $h_{p,s}:=p^{\lceil\log_ps\rceil}\quad(s\ge1)$.
Fix a shift vector $a \ne 0$. 
The cosets of $\langle a\rangle=\{\lambda a:\lambda\in\mathbb F_p\}$ partition $V_{p,n}$ into $M$ affine lines, where on each line the edges $x\to x+a$ form a directed $p$-cycle.
We call a point as \textit{marked} if it belongs to the set $B$ and \textit{unmarked} otherwise.

Suppose that a cycle has marked and unmarked points.
We will call a maximal sequence of consecutive marked points as \textit{run}.
Thus each run starts just after an unmarked point and ends just before an unmarked point.
For example, in a $5$-cycle, if positions numbered $0,1,3$ are marked, the runs are $0,1$ and $3$.
The first run $\{0,1\}$ contributes one edge, and $\{3\}$ contributes no edge. 
If positions $0,1,4$ are marked instead, then they form one run. 
In general, a run with $\ell$ points contributes $\ell-1$ edges with both endpoints in $B$. 
Hence a cycle with $k$ marked points and $j$ runs contributes $k-j$ to $N_a(B)$.

With these terms, we are now ready to state and prove our theorem.
\begin{samepage}
\begin{theorem}\label{thm:odd-prime}
Let $p$ be an odd prime and $n\ge1$.
For $0\le t\le N_p$, we have the following.
\begin{enumerate}
\item If $0\le s\le2M$, then $f_{p,n}(s,t)=s\,c_{p,n}(t)$.
\item If $1\le s\le N_p$, then
$g_{p,n}(s,t)=st$ if and only if $h_{p,s}\mid t$.
\item If $N_p-2M\le s\le N_p$, then
$g_{p,n}(s,t)=t^2-(N_p-s)c_{p,n}(t)$.
\end{enumerate}
\end{theorem}
\end{samepage}

\begin{proof}
Fix a shift vector $a \ne 0$. 
Consider a $p$-cycle and suppose that $0<k<p$ of its points belong to $B$.
Let $j$ be the number of runs. Each run contains at least one
marked point and is followed by an unmarked point, so
$j\le k$ and $j\le p-k$. Thus the contribution to $N_a(B)$ is
\[
 k-j\ge k-\min\{k,p-k\}=\max\{0,2k-p\}.
\]

If $k\le p-k$, we can choose the marked points so that no two
are consecutive. Then there are $k$ runs, giving overlap zero.
If $k>p-k$, we can choose the unmarked points so that no two
are consecutive. The marked points between successive unmarked
points then form $p-k$ runs, giving overlap $2k-p$.
An empty cycle contributes zero, and a full cycle contributes $p$.
Therefore the minimum overlap on a cycle with $k$ marked points is
\[
 \phi_p(k)=\max\{0,2k-p\}.
\]

A cycle has overlap zero exactly when no two marked points are consecutive. 
Each marked point must then be followed by an unmarked point, so $k\le p-k$, or $k\le q$.
We can attain zero overlap by marking any $k\le q$ of the positions $0,2,\ldots,2q-2$. 
Thus each cycle can contain at most $q$ marked points with zero overlap.

We already have that
\[
 \phi_p(k)\ge0,\qquad \phi_p(k)\ge2k-p.
\]
We also have $\phi_p(k)\ge k-q$. This holds when $k\le q$ because
$k-q\le0$. When $k\ge q+1$, we have
\[
 \phi_p(k)=2k-p=(k-q)+(k-q-1)\ge k-q.
\]
Let $k_L=|B\cap L|$ for each $p$-cycle $L$. 
Since there are $M$ number of cycles
and $\sum_L k_L=t$, summing these three bounds gives
\[
 N_a(B)\ge\max\{0,t-qM,2t-pM\}=c_{p,n}(t).
\]
We now construct the set $B$ that attains the lower bound for $2M$ shifts. 
Let $V_{p,n}=U\times\F_p$, where $U=\F_p^{n-1}$ and $|U|=M$. 
Order the layers as follows.
\[
 0,2,4,\ldots,2q,1,3,\ldots,2q-1.
\]
A shift by the element $(u,1)$ translates layer $i$ to layer $i+1$, with the indices in modulo $p$.

For $t<N_p$, let $t=mM+r$ with $0\le r<M$.
Fill the first $m$ layers and place an arbitrary $r$-point set in
the next layer. At $t=N_p$, fill every layer.
For every $u\in U$, the overlap when the shift is $(u,1)$ is
\[
 \begin{cases}
 0,&m<q,\\
 r,&m=q,\\
 (2m-p)M+2r,&q+1\le m\le p-1.
 \end{cases}
\]
Shifting by $(u,1)$ sends each of the first $q$ layers outside those layers. 
The next layer, indexed by $2q$, shifts into layer $0$, so choosing $r$ points in it contributes $r$ to the overlap.

For each later layer $i$, both layers $i-1$ and $i+1$ have already been included. 
Choosing $r$ points in layer $i$ gives overlap $r$ from layer $i-1$ into layer $i$, and another $r$ from layer $i$ into layer $i+1$. 
These counts do not depend on $u$ or on which $r$ points we choose.

At $t=N_p$, we have $B=V_{p,n}$, so the overlap is $N_p$.
Thus the overlap equals $c_{p,n}(t)$ in every case.
Since $N_a(B)=N_{-a}(B)$, the same holds for shifts $(u,-1)$.
We therefore have $2M$ shifts attaining the one-shift minimum.
The zero shift also satisfies the lower bound.
Choosing $A$ to be any $s$ of these $2M$ shifts gives
$r(A,B,B)=s\,c_{p,n}(t)$ and proves the first part of the theorem.

Now for the second statement.
$\Stab(B)$ is an $\F_p$-subspace, and
$r(A,B,B)=st$ holds exactly when $A\subseteq\Stab(B)$.
$B$ is a union of cosets of $\Stab(B)$, so $|\Stab(B)|$ divides $t$.
The rest of the argument is the same as in \Cref{thm:upper-main}.

Finally, \Cref{prop:duality} holds in any finite abelian group.
If $s\ge N_p-2M$, then $N_p-s\le2M$.
Replace $A$ by $V_{p,n}\setminus A$ and apply the first part of the theorem.
\end{proof}

For $p=3$, the lower endpoint formula holds for $s\le2N/3$,
and the upper endpoint formula holds for $s\ge N/3$.
For larger $p$, the corresponding parameter ranges are
$s\le2N_p/p$ and $s\ge(p-2)N_p/p$, respectively.

\section{Further discussion}\label{sec:discussion}

We have determined the lower endpoint $f_n(s,t)$ for
$s\le2N/3$ and the upper endpoint $g_n(s,t)$ for $s\ge N/3$.
As an application of these endpoint results, we have also improved the range of Schur-triple minimum. 
We now summarize the endpoint results. 
In the next subsection we will provide
constructions and holes inside the spectrum.

\subsection{The remaining endpoint problem}

Let $N=3^n$, $0\le s,t\le N$, and $c_n(t)=\max\{0,t-N/3,2t-N\}$.
For $u\ge1$, let $h_u=3^{\lceil\log_3u\rceil}$.
Also, let $\varepsilon(t)=0$ when $3\mid t$ and
$\varepsilon(t)=1$ otherwise.

\begin{table}[htbp]
\centering
\small
\renewcommand{\arraystretch}{1.25}
\begin{tabular}{@{}lp{0.42\textwidth}p{0.40\textwidth}@{}}
\toprule
Endpoint & Range & Formula \\
\midrule
$f_n(s,t)$
& $0\le s\le2N/3$
& $s\,c_n(t)$ \\

$g_n(s,t)$
& $N/3\le s\le N$
& $t^2-(N-s)c_n(t)$ \\

$g_n(s,t)$
& $1\le s\le N$, $h_s\mid t$
& $st$ \\

$g_n(s,t)$
& $1\le t\le N$, $h_t\le s\le N$
& $t^2$ \\

$g_n(s,t)$
& $1\le s\le3$
& $st-(s-1)\varepsilon(t)$ \\

$g_n(s,t)$
& $4\le s<N/3$, $h_s\nmid t$
& \Cref{thm:remaining-endpoint-bounds} \\
\bottomrule
\end{tabular}
\caption{Endpoint formulas and bounds in $\F_3^n$.}
\label{tab:endpoint-summary}
\end{table}

The identity $f_n(N-s,t)=t^2-g_n(s,t)$ relates the lower endpoint with the upper endpoint.
It therefore remains to determine $g_n(s,t)$ in the unresolved range.
In \Cref{rem:remaining-bounds-sharp} and \Cref{rem:distinct-bounds-exact}, the lower bound is exact. 
This leads us to ask whether this always holds.

\begin{problem}\label{prob:endpoints}
Let $4\le s<N/3$, let $h=h_s=3^k$, and suppose that
$0<t<N$ and $h\nmid t$.
Let $t=mh+r$, where $0<r<h$, and let
\[
 c_k(r)=\max\{0,r-h/3,2r-h\}.
\]
Is it true that
\[
 g_n(s,t)
 =s(t-r)+r^2-(h-s)c_k(r)?
\]
\end{problem}

\subsection{Values inside the endpoints}\label{sec:middle}

Knowing the endpoint values, does not determine which values occur between them.
An interesting direction is to determine which additive
structures produce intervals, arithmetic progressions, and
holes in the spectrum.

For a fixed set $B\subseteq V_n$, let
\[
 \mathcal R_s(B)
 :=
 \{r(A,B,B):A\subseteq V_n,\ |A|=s\}.
\]
Then
\[
 r_n(s,t)=\bigcup_{\substack{B\subseteq V_n\\|B|=t}}
 \mathcal R_s(B).
\]
The following constructions describe $\mathcal R_s(B)$
for two different choices of $B$.

\subsubsection{Subspaces}

If $L\le V_n$ is a subspace and has size $t=3^d$, then
\[
 N_a(L)=
 \begin{cases}
 t,&a\in L,\\
 0,&a\notin L.
 \end{cases}
\]

\begin{proposition}\label{prop:subspace-values}
Let $L\le V_n$ have size $t=3^d$, where $0\le d\le n$, and let $j \in \mathbb{Z}$.
Then
\[
 \mathcal R_s(L)
 =
 \left\{
 jt:
 \max(0,s-(N-t))\le j\le\min(s,t)
 \right\}.
\]
In particular, all these values belong to $r_n(s,t)$.
\end{proposition}

\begin{proof}
For $B=L$, we have $r(A,L,L)=t|A\cap L|$.
A set $A$ of size $s$ can contain exactly $j$ elements of $L$, when $0\le j\le t$ and $0\le s-j\le N-t$.
\end{proof}

\subsubsection{Sidon sets}

We will now use Sidon sets to give intervals of attainable values.
A set $B\subseteq V_n$ is \emph{Sidon} if each nonzero ordered
difference occurs at most once.
Equivalently, for all $b_1,b_2,b_3,b_4\in B$,
\[
 b_1-b_2=b_3-b_4\ne0
 \quad\Longrightarrow\quad
 b_1=b_3\text{ and }b_2=b_4.
\]

If $B\subseteq\F_3^n$ is Sidon, then it has $|B|(|B|-1)$ nonzero ordered differences, and they are all distinct. Hence $|B|(|B|-1)\le 3^n-1$.
This gives an upper bound on the size of such a set.
For a nice overview of constructions of Sidon sets see~\cite{EM} and the references therein.
Since every subset of a Sidon set is a Sidon set, these constructions also provide examples of every smaller cardinality.

For a Sidon set $B$ of size $t$, let $D:=\{b-b':b,b'\in B,\ b\ne b'\}$.
Then
\[
 |D|=t(t-1),
 \qquad
 N_a(B)=
 \begin{cases}
 t,&a=0,\\
 1,&a\in D,\\
 0,&a\notin D\cup\{0\}.
 \end{cases}
\]
Therefore choosing $A$ results in selecting shifts of weights
$t$, $1$, and $0$.

\begin{proposition}\label{prop:sidon-values}
Suppose that $B\subseteq V_n$ is Sidon with $|B|=t$.
Let $d=t(t-1)$,  $z=(N-1-d)$ and $k \in \mathbb{Z}$
Then
\begin{align*}
 \mathcal R_s(B)
 ={}&
 \{k:\max(0,s-z)\le k\le\min(s,d)\}\\
 &{}\cup
 \{t+k:\max(0,s-1-z)\le k\le\min(s-1,d)\}.
\end{align*}
\end{proposition}

\begin{proof}
If $0\notin A$, choose $k$ shifts from $D$ and $s-k$
from the $z$ zero-overlap shifts.
This is possible exactly when $0\le k\le d$ and
$0\le s-k\le z$, and gives count $k$.

If $0\in A$, the zero shift contributes $t$.
Apply the same argument to the remaining $s-1$ shifts.
\end{proof}

In particular, if $1\le s\le\min(d,z)$, then
\begin{equation}\label{eq:sidon-intervals}
 \{0,1,\ldots,s\}\cup\{t,t+1,\ldots,t+s-1\}
 \subseteq r_n(s,t).
\end{equation}
Note that a gap in $\mathcal R_s(B)$ for this particular $B$ is not
necessarily a hole in $r_n(s,t)$. 

\subsection{Consecutive holes}

In the previous subsection we discussed constructions to attain values in the spectrum. 
We now show that certain values cannot occur for any choice of $B$.

\begin{proposition}\label{prop:near-top-hole}
If $3\mid t$, then $t-1\notin r_n(1,t)$.
If $3\le t\le N-3$ and $3\mid t$, then $t-1$ is a hole
strictly between the endpoints.
\end{proposition}

\begin{proof}
The zero shift contributes $t$.
For $a\ne0$, let $u_i$ count the cosets of $\langle a\rangle$
containing exactly $i$ points of $B$.
Using \eqref{eq:cycle-rule}, we obtain
\[
 t=u_1+2u_2+3u_3,
 \qquad N_a(B)=u_2+3u_3.
\]
Subtracting the two equations gives $t-N_a(B)=u_1+u_2$. 
If $N_a(B)=t-1$, then $u_1+u_2=1$. This gives $t\equiv1$ or $2\pmod3$, contrary to $3\mid t$.

If $3\le t\le N-3$ and $3\mid t$, then
\[
 f_n(1,t)=c_n(t)\le t-3<t-1<t=g_n(1,t),
\]
so the excluded value is a hole.
\end{proof}

The following proposition tells us that a block of holes appear below the maximum.

\begin{proposition}\label{prop:top-blocks}
Let $s\ge3$ be a power of $3$, and suppose that
$s\mid t$ and $0<t<N$.
Then $g_n(s,t)=st$, and every value in
\[
 \{st-s+1,\ldots,st-1\}
\]
is a hole in $r_n(s,t)$.
\end{proposition}

\begin{proof}
Fix $A,B$ of sizes $s,t$, and let $H=\Stab(B)$ and $\ell=|H|$.
Suppose that $d=st-r(A,B,B)>0$. 
We show that $d\ge s$.
Since both $B$ and $B+a$ are unions of $H$-cosets, their intersection is also a union of $H$-cosets.
Therefore $N_a(B)$ is divisible by $\ell$ for every $a$.
We call $t-N_a(B)$ the \emph{deficit} of the shift $a$.
The total deficit is
\[
 d=st-r(A,B,B)
   =\sum_{a\in A}(t-N_a(B)).
\]
Since $\ell\mid t$, every deficit of a shift and the total deficit
are divisible by $\ell$.
\begin{enumerate}
    \item If $\ell\ge s$, then $t-N_a(B)>0 \implies t-N_a(B)\ge\ell \implies t-N_a(B)\ge s$.
    \item If $\ell=1$, then $H=\{0\}$. Thus for $a \ne 0$ we have $B+a\ne B$. 
        And so their intersection has size less than $t$. 
        \Cref{prop:near-top-hole} excludes $N_a(B)=t-1$, and so deficit of at least $2$.
        Hence, $d\ge2(s-1)\ge s$.
    \item If $1<\ell<s$, then $3\le\ell\le s/3$.
At most $\ell$ selected shifts stabilize $B$, and every
other shift has deficit at least $\ell$.
Hence, $d\ge\ell(s-\ell)\ge2s$.
\end{enumerate}

Thus it is not possible to have $1,\ldots,s-1$ as deficit.

For attaining $st$, choose a subspace $K$ of size $s$, take $B$ to be a union of $t/s$ cosets of $K$, and set $A=K$.
Until now we have shown that $st$ is attainable and $\{st-s+1,\ldots,st-1\}$ is not. 
We will now show an attainable value below these excluded values.
Choose $a\notin\Stab(B)$ and replace one element of $A$ by $a$.
The overlap $N_a(B)$ is divisible by $s$ and smaller than $t$,
so the new count is at most $st-s$.
The excluded values therefore lie between
attainable values.
\end{proof}

Taking $t=s$ and $s<N$ shows that these blocks have
arbitrarily large length as the dimension grows.

The construction using a subspace suggests that unions of cosets, possibly of subspaces of different dimensions, form the natural family to investigate next.
We end our paper with the following problem.

\begin{problem}\label{prob:coset-spectrum}
Let $1\le s,t\le N$ and $k\ge1$.
Which values in $r_n(s,t)$ can be attained by pairs $(A,B)$
for which $B$ is a disjoint union of at most $k$ affine subspaces
of $V_n$, possibly of different dimensions? 

\end{problem}
\medskip

\textbf{Acknowledgments.}
The second author thanks Jonathan Jedwab for his inspiring talk on \textit{Additive Triples in Groups of Odd Prime Order} at Combinatorics 2024 in Carovigno, which led to the collaboration between the two authors.
\medskip

\textbf{AI declaration.} 
The second author used AI for literature searches and for assistance with exposition.

\end{document}